\documentclass[10pt,onecolumn]{article}

\usepackage{times}
\usepackage[dvips]{graphicx}
\usepackage[hang,scriptsize]{subfigure}
\usepackage{latexsym}

\usepackage{amsmath} 
\usepackage{amssymb} 
\usepackage{graphics} 
\usepackage{listings}    
\usepackage{stmaryrd}    

\makeatletter
\def\@maketitle{%
  \vbox to 6.5cm{%
    \hsize\textwidth
    \linewidth\hsize
    \vspace{1.5cm}
    \centering
    {\bfseries\LARGE \@title \par}
    \vspace{12pt}
    {\fontsize{11pt}{13pt}\selectfont \begin{tabular}[t]{c}\@author \end{tabular}\par}
    \vfill} 
}
\renewcommand\section{\@startsection{section}{1}{\z@}%
                       {-12\p@ \@plus -4\p@ \@minus -4\p@}%
                       {6\p@ \@plus 4\p@ \@minus 4\p@}%
                       {\normalfont\large\bfseries
                        \rightskip=\z@ \@plus 8em\pretolerance=10000 }}
\renewcommand\subsection{\@startsection{subsection}{2}{\z@}%
                       {-12\p@ \@plus -4\p@ \@minus -4\p@}%
                       {6\p@ \@plus 4\p@ \@minus 4\p@}%
                       {\normalfont\fontsize{11pt}{13pt}\selectfont\bfseries
                        \rightskip=\z@ \@plus 8em\pretolerance=10000 }}
\renewcommand\subsubsection{\@startsection{subsubsection}{3}{\z@}%
                       {-12\p@ \@plus -4\p@ \@minus -4\p@}%
                       {6\p@ \@plus 4\p@ \@minus 4\p@}%
                       {\normalfont\normalsize\itshape}}
\renewcommand\paragraph{\@startsection{paragraph}{4}{\z@}%
                       {-12\p@ \@plus -4\p@ \@minus -4\p@}%
                       {-0.5em \@plus -0.22em \@minus -0.1em}%
                       {\normalfont\normalsize\itshape}}
\makeatother

\renewenvironment{abstract}%
  {\small
    \list{}{\labelwidth0pt
      \leftmargin0pt \rightmargin\leftmargin
      \listparindent\parindent \itemindent0pt
      \parsep0pt
      }%
    \item[\hskip\labelsep\bfseries\abstractname\enspace --] \itshape}{\endlist}

\newcommand{\keywordsname}{Keywords}
\newenvironment{keywords}%
  {\small
    \list{}{\labelwidth0pt
      \leftmargin0pt \rightmargin\leftmargin
      \listparindent\parindent \itemindent0pt
      \parsep0pt
      }%
    \item[\hskip\labelsep\bfseries\keywordsname:]}{\endlist}

\begin{document}

\title{Fusion of imprecise, uncertain and conflicting beliefs\\
 with DSm rules of combination}

\author{\begin{tabular}{c@{\extracolsep{8em}}c}
{\bf Jean Dezert} & {\bf Florentin Smarandache}\\
ONERA & Department of Mathematics\\
9 Av. de la  Division Leclerc & University of New Mexico\\
92320 Ch\^{a}tillon & Gallup, NM 8730\\
France & U.S.A.\\
{\tt Jean.Dezert@onera.fr} & {\tt smarand@unm.edu}
\end{tabular}}

\date{}
\maketitle
\pagestyle{plain}

\begin{abstract}
In this paper one studies, within Dezert-Smarandache Theory (DSmT), the case when the sources of information provide imprecise belief functions/masses, and we generalize the DSm rules of combination (classic or hybrid rules) from scalar fusion to sub-unitary interval fusion and, more general, to any set of sub-unitary interval fusion. This work generalizes previous works available in literature which appear limited to IBS (Interval-valued belief structures) in the Transferable Belief Model framework.  Numerical didactic examples of these new DSm fusion rules for dealing with imprecise information are also presented.
\end{abstract}



\section{Introduction}

During the last two years, we have developed a new mathematical theory, the DSmT (Dezert Smarandache Theory),  for combining uncertain and conflicting sources of information  \cite{Dezert_2002b,Dezert_2003,Dezert_2004Book,Smarandache_2002}. The DSmT is based on a new modeling of the fusion problem and propose new rules of combination which appear to be more attractive than the classical Dempster's rule of combination proposed by G. Shafer within the development of the Dempster-Shafer Theory (DST) \cite{Shafer_1976}, specially when one has to deal with high conflicting sources of information and/or dynamical fusion problems, where the frame of discernment changes with time.  The DSmT allows  the fusion of sources, thanks to the classical DSm rule, on free-DSm models (model where all hypotheses of the frame $\Theta$ are partially overlapping without possibility for refinement), but more generally on any more complex/restricted model (like the Shafer's model) including any kind of integrity constraints thanks to the DSm hybrid rule of combination.\\

 Until now,  we had focused our efforts on the fusion of {\it{precise}} uncertain and conflicting/paradoxist generalized basic belief assignments (gbba). We mean here by precise gbba, basic belief functions/masses $m(.)$ defined precisely on hyper-power set  $D^\Theta$ where each mass $m(X)$, where $X$ belongs to $D^\Theta$, is represented by only one real number belonging to $[0,1]$ such that $\sum_{X\in D^\Theta}m(X)=1$. In this paper, we extend the DSm fusion rules for dealing with {\it{admissible imprecise generalized basic belief assignments}} $m^I(.)$ defined as real subunitary intervals of $[0,1]$, or even more general as real subunitary sets [i.e. 
sets, not necessarily intervals]. An imprecise belief assignment $m^I(.)$ over $D^\Theta$ is said admissible if and only if there exists for every $X\in D^\Theta$ at least one real number $m(X)\in m^I(X)$ such that $\sum_{X\in D^\Theta}m(X)=1$. The idea to work with imprecise belief structures represented by real subset intervals of $[0,1]$ is not new and we strongly encourage the reader to examine previous Lamata \& Moral's together with Den\oe ux's works for instance on this topic in \cite{Lamata_1994,Denoeux_1997,Denoeux_1999} and references therein. The proposed works available in the literature, upon our knowledge were limited only to sub-unitary interval combination in the framework of Transferable Belief Model (TBM) developed by Smets \cite{Smets_1994,Smets_2000}. We extend Lamata \& Moral's together with Den\oe ux's subunitary interval-valued masses to subunitary set-valued masses; therefore the 
closed intervals used by Denoeux to denote imprecise masses are generalized to any sets included in [0, 1], i.e. in our case these sets can be unions of (closed, open, or half-open / half-closed) intervals and/or scalars all in $[0,1]$. In this work, the proposed extension is done in the context of DSmT framework, although it can also apply directly to fusion of IBS within TBM as well if the user prefers to adopt TBM rather than DSmT.\\

In many fusion problems, it seems very difficult (if not impossible) to have precise sources of evidence generating precise basic belief assignments (specially when belief functions are provided by human experts), and a more flexible plausible and paradoxical  theory supporting imprecise information becomes necessary. This paper proposes a new issue to deal with the fusion of imprecise, uncertain and conflicting source of information. The section \ref{Sec:CPB} presents briefly the DSm rule of combination for precise belief functions. In section \ref{Sec:SetOperators}, we present the operations on 
sets for the paper to be self-contained and necessary to deal with imprecise nature of information in our framework. In section \ref{Sec:CIB}, we propose an issue to combine simple imprecise belief assignment corresponding only to sub-unitary intervals also known as IBS (Interval-valued belief structures) in \cite{Denoeux_1997}. In section \ref{Sec:GIB}, we present the generalization of our new fusion rules to combine any type of imprecise belief assignment which may be represented by the union of several sub-unitary (half-) open intervals, (half-)closed intervals and/or sets of points belonging to [0,1]. Several numerical examples are also given. In the sequel, one uses the notation $(a,b)$ for an open interval, $[a,b]$ for a closed interval, and $(a,b]$ or $[a,b)$ for a half open and half closed interval.

\section{Combination of precise beliefs}
\label{Sec:CPB}

\subsection{General DSm rule of combination}

Let's consider a frame of discernment of a fusion problem $\Theta=\{\theta_1,\theta_2,\ldots,\theta_n\}$, its hyper-power set $D^\Theta$ (i.e. the set of all propositions built from elements $\theta_i$ of $\Theta$ with $\cap$ and $\cup$ operators \cite{Dezert_2003f,Dezert_2004Book}, and $k$ independent (precise) sources of information $\mathcal{B}_1$, $\mathcal{B}_2$, $\ldots$, $\mathcal{B}_k$ with their associated generalized basic belief assignments (gbba)  
$m_1(.)$, $m_2(.)$, $\ldots$, $m_k(.)$ defined over  $D^\Theta$. Let $\mathbf{M}$ be the mass matrix
\begin{equation*}
\mathbf{M}=
\begin{bmatrix}
m_{11} & m_{12} & \ldots & m_{1d}\\
m_{21} & m_{22} & \ldots & m_{2d}\\
 \ldots &  \ldots & \ldots & \ldots\\
m_{k1} & m_{k2} & \ldots & m_{kd}\\
\end{bmatrix}
\end{equation*}
\noindent where $d={\mid D^\Theta\mid}$ is the dimension of the hyper-power set, and $m_{ij}\in [0,1]$ for all $1\leq i\leq k$ and $1\leq  j\leq d$, is the mass assigned by source $\mathcal{B}_i$ to the element $A_j\in D^\Theta$. We use the DSm ordering procedure presented in \cite{Dezert_Smarandache_2003,Dezert_2004Book}
for enumerating the elements $A_1$, $A_2$, \ldots, $A_d$ of the hyper-power set $D^\Theta$.
The matrix $\mathbf{M}$ characterizes all information available which has to be combined to solve the fusion problem under consideration. Since $m_1(.)$, $m_2(.)$, $\ldots$, $m_k(.)$ are gbba, the summation on each raw of the matrix must be one.  For any (possibly hybrid) model ${\mathcal{M}(\Theta)}$, we apply the DSm general rule of combination (also called DSm hybrid rule) for $k\geq 2$ sources to fusion the masses \cite{Dezert_2004Book}  
defined for all $A\in D^\Theta$ as:
\begin{equation}
m_{\mathcal{M}(\Theta)}(A)\triangleq 
\phi(A)\Bigl[ S_1(A) + S_2(A) + S_3(A)\Bigr]
 \label{eq:DSmHkno1}
\end{equation}
\noindent
$\phi(A)$ is the {\it{characteristic emptiness function}} of the set $A$, i.e. $\phi(A)= 1$ if  $A\notin \boldsymbol{\emptyset}$ and $\phi(A)= 0$ otherwise. $\boldsymbol{\emptyset}\triangleq \{\emptyset,\boldsymbol{\emptyset}_{\mathcal{M}}\}$ represents the set absolutely empty and of all relatively empty  elements belonging to $D^\Theta$ (elements/propositions which have been forced to empty set in the chosen hybrid model ${\mathcal{M}(\Theta)}$). If no constraint is introduced in the model, $\boldsymbol{\emptyset}$ reduces to $\{\emptyset\}$ and this corresponds to the free-DSm model  \cite{Dezert_2004Book}.
If all constraints of exclusivity between elements $\theta_i\in\Theta$ are introduced, the hybrid model  ${\mathcal{M}(\Theta)}$ corresponds to the Shafer's model on which is based the Dempster-Shafer Theory (DST) \cite{Shafer_1976}. $S_1(A)$, $S_2(A)$ and $S_3(A)$ are defined by

\begin{equation}
S_1(A)\triangleq \sum_{\overset{X_1,X_2,\ldots,X_k\in D^\Theta}{(X_1\cap X_2\cap\ldots\cap X_k)=A}} \prod_{i=1}^{k} m_i(X_i)
\label{eq:S1}
\end{equation}

\begin{equation}
S_2(A)\triangleq \sum_{\overset{X_1,X_2,\ldots,X_k\in\boldsymbol{\emptyset}}{ [\mathcal{U}=A]\vee [(\mathcal{U}\in\boldsymbol{\emptyset}) \wedge (A=I_t)]}} \prod_{i=1}^{k} m_i(X_i)\end{equation}

\begin{equation}
S_3(A)\triangleq\sum_{\overset{X_1,X_2,\ldots,X_k\in D^\Theta}{\overset{(X_1\cup X_2\cup\ldots\cup X_k)=A}{\overset{(X_1\cap X_2\cap \ldots\cap X_k)\in\boldsymbol{\emptyset}}{}}}}  \prod_{i=1}^{k} m_i(X_i)
\end{equation}

where $I_t\triangleq\theta_1\cup \theta_2\cup \ldots\cup \theta_n$ and  $ \mathcal{U}\triangleq u(X_1)\cup u(X_2)\cup \ldots \cup u(X_k)$. $u(X)$ is the union of all singletons $\theta_i$ that compose $X$.  For example, if $X$ is a singleton then $u(X)=X$; if  $X=\theta_1\cap \theta_2$ or $X=\theta_1\cup \theta_2$ then $u(X)=\theta_1\cup \theta_2$; if $X=(\theta_1\cap \theta_2)\cup \theta_3$ then $u(X)=\theta_1\cup \theta_2\cup\theta_3$, etc; by convention $u(\emptyset)\triangleq\emptyset$.

\clearpage
\newpage

\subsection{Examples}
\label{subsec:Examples2.2}
Let's consider at time $t$ the frame of discernment $\Theta=\{ \theta_1,\theta_2,\theta_3\}$ and  two independent bodies of evidence $\mathcal{B}_1$ and $\mathcal{B}_2$ with the generalized basic belief assignments $m_1(.)$ and $m_2(.)$ given by:

\begin{table}[h]
\begin{equation*}
\begin{array}{|c|c|c|}
\hline
A\in D^\Theta & m_1(A) & m_2(A) \\
\hline
\theta_1 & 0.1 & 0.5\\
\theta_2 & 0.2 & 0.3\\
\theta_3 & 0.3 & 0.1\\
\theta_1\cap\theta_2 & 0.4 & 0.1\\
\hline
\end{array}
\end{equation*}
\caption{Inputs of the fusion with precise bba}
\label{mytable1}
\end{table}

Based on the free DSm model and the classical DSm rule \eqref{eq:S1}, the combination denoted by the symbol $\oplus$ (i.e. $m(.)=[m_1\oplus m_2](.)$) of these two precise sources of evidence is
\begin{table}[h]
\begin{equation*}
\begin{array}{|c|c|}
\hline
A\in D^\Theta & m(A)=[m_1\oplus m_2](A) \\
\hline
\theta_1 & 0.05\\
\theta_2 & 0.06\\
\theta_3 & 0.03\\
\theta_1\cap\theta_2 & 0.52\\
\theta_1\cap\theta_3 & 0.16\\
\theta_2\cap\theta_3 & 0.11\\
\theta_1\cap\theta_2 \cap \theta_3& 0.07\\
\hline
\end{array}
\end{equation*}
\caption{Fusion with DSm classic rule}
\label{mytable2}
\end{table}
Then, assume at time $t+1$ one finds out for some reason that the free-DSm model has to be changed by introducing the constraint $\theta_1\cap\theta_2=\emptyset$ which involves also $\theta_1\cap\theta_2\cap\theta_3=\emptyset$. This characterizes the hybrid-model $\mathcal{M}$ we have to work with. Then one uses the general DSm hybrid rule of combination for scalars (i.e. for precise masses $m_1(.)$ and $m_2(.)$ to get the new result of the fusion at time $t+1$. According to \eqref{eq:DSmHkno1}, one obtains $m(\theta_1\cap\theta_2\overset{\mathcal{M}}{\equiv}\emptyset)=0$, $m(\theta_1\cap\theta_2 \cap \theta_3\overset{\mathcal{M}}{\equiv}\emptyset )=0$ and

\begin{table}[h]
\begin{equation*}
\begin{array}{|c|c|}
\hline
A\in D^\Theta & m(A)\\
\hline
\theta_1 & 0.05+[0.1(0.1)+0.5(0.4)]=0.26\\
\theta_2 & 0.06+[0.2(0.1)+0.3(0.4)]=0.20\\
\theta_3 & 0.03+[0.3(0.1)+0.1(0.4)]=0.10\\
\theta_1\cap\theta_3 & 0.16\\
\theta_2\cap\theta_3 & 0.11\\
\theta_1\cup\theta_2 & 0 + [0.13] + [0.04]=0.17\\
\hline
\end{array}
\end{equation*}
\caption{Fusion with DSm hybrid rule for model $\mathcal{M}$}
\label{mytable3}
\end{table}

\section{Operations on sets}
\label{Sec:SetOperators}

To manipulate imprecise information and for the paper to be self-contained, we need to introduce operations on sets as follows (detailed presentations on Interval Analysis and Methods can be found in \cite{Hayes_2003,Jaulin_2001,Moore_1966,Moore_1979,Neumaier_1990}). The interval operations defined here about imprecision are similar to the rational interval extension through the interval arithmetics \cite{SIGLA}, but they are different from Modal Interval Analysis which doesn't serve our fusion needs.  We are not interested in a dual of an interval $[a, b]$, used in the Modal Interval Analysis, because we always consider $a \leq b$, while its dual, $\text{Du}([a, b])=[b, a]$, doesn't occur. Yet, we generalize the interval operations to any set operations.  Of course, for the fusion we only need real sub-unitary sets, but these defined set operations can be used for any kind of sets.\\

Let $S_1$ and $S_2$ be two (unidimensional) real standard subsets of the unit interval $[0, 1]$, and a number $k\in [0,1]$, then one defines \cite{Smarandache_2002} :

\begin{itemize}
 \item
 {\bf{Addition of sets}}
 \begin{equation*}
 S_{1}\boxplus S_{2} =S_{2}\boxplus S_{1}\triangleq \{ x \mid x = s_{1}+s_{2},  s_{1} \in 
S_{1},s_{2} \in S_{2} \}
 \quad\text{with} \quad 
 \begin{cases}
\inf(S_{1}\boxplus S_{2})=\inf(S_1) + \inf(S_2)\\
\sup(S_{1}\boxplus S_{2})=\sup (S_1) + \sup(S_2)
\end{cases}
\label{eq:addition}
\end{equation*}

\noindent
and, as a particular case, we have 
 \begin{equation*}
 \{k\}\boxplus S_2=S_2\boxplus\{k\}=\{ x \mid x = k +s_{2}, s_{2} \in S_{2} \}
 \quad\text{with} \quad 
\begin{cases}
\inf(\{k\}\boxplus S_{2})= k + \inf(S_2)\\
\sup(\{k\}\boxplus S_{2}) = k + \sup(S_2)
\end{cases}
\end{equation*}

{\it{Examples}}: 

$[0.1,0.3]\boxplus[0.2,0.5]=[0.3,0.8]$ because $0.1+0.2=0.3$ and $0.3+0.5=0.8$; 

$(0.1,0.3]\boxplus[0.2,0.5]=(0.3,0.8]$; 

$[0.1,0.3]\boxplus(0.2,0.5]=(0.3,0.8]$; 

$[0.1,0.3)\boxplus[0.2,0.5]=[0.3,0.8)$; 

$[0.1,0.3]\boxplus[0.2,0.5)=[0.3,0.8)$; 

$(0.1,0.3]\boxplus(0.2,0.5)=(0.3,0.8)$; 

$[0.7,0.8]\boxplus[0.5,0.9]=[1.2,1.7]$; 

$\{0.4\}\boxplus[0.2,0.5]= [0.2,0.5]\boxplus\{0.4\}=[0.6,0.9]$ because $0.4+0.2=0.6$ and $0.4+0.5=0.9$; 

$\{0.4\}\boxplus(0.2,0.5]=(0.6,0.9]$;

$\{0.4\}\boxplus[0.2,0.5)=[0.6,0.9)$;

$\{0.4\}\boxplus(0.2,0.5)=(0.6,0.9)$.

\item
{\bf{Subtraction of sets}}
 \begin{equation*}
 S_{1}\boxminus S_{2} \triangleq \{ x \mid x = s_{1}-s_{2},  s_{1} \in 
S_{1}, s_{2} \in S_{2} \}
 \quad\text{with} \quad 
\begin{cases}\inf(S_{1}\boxminus S_{2})=\inf(S_1) - \sup(S_2)\\
\sup(S_{1}\boxminus S_{2})=\sup(S_1) - \inf(S_2)
\end{cases}
\label{eq:addition}
 \end{equation*}
\noindent
and, as a particular case, we have 
 \begin{equation*}
\{k\}\boxminus S_2=\{ x \mid x = k -s_{2},s_{2} \in S_{2} \}
 \quad\text{with} \quad 
\begin{cases}
\inf(\{k\}\boxminus S_{2})= k - \sup(S_2)\\
\sup(\{k\}\boxminus S_{2}) = k - \inf(S_2)\end{cases}
 \end{equation*}
\noindent
and similarly for $S_2\boxminus\{k\}$ with
$\begin{cases}
\inf(S_2\boxminus\{k\})= \inf(S_2) - k\\
\sup(S_2\boxminus\{k\}) = \sup(S_2) - k
\end{cases}$\\

{\it{Examples}}: 

$[0.3,0.7]\boxminus[0.2,0.3]=[0.0,0.5]$ because $0.3-0.3=0.0$ and $0.7-0.2=0.5$; 

$[0.3,0.7]\boxminus\{0.1\}=[0.2,0.6]$; 

$\{0.8\}\boxminus[0.3,0.7]=[0.1,0.5]$ because $0.8-0.7=0.1$ and $0.8-0.3=0.5$; 

$[0.1,0.8]\boxminus[0.5,0.6]=[-0.5,0.3]$; 

$[0.1,0.8]\boxminus[0.2,0.9]=[-0.8, 0.6]$; 

$[0.2,0.5]\boxminus[0.1,0.6]=[-0.4, 0.4]$.

  \item
 {\bf{Multiplication of sets}}
 \begin{equation*}
S_{1}\boxdot S_{2} \triangleq \{ x \mid x = s_{1}\cdot s_{2}, s_{1} \in 
S_{1},s_{2} \in S_{2} \}
 \quad\text{with} \quad 
\begin{cases}
\inf(S_{1}\boxdot S_{2})=\inf(S_{1})\cdot \inf(S_{2})\\
\sup(S_{1}\boxdot S_{2})=\sup(S_{1})\cdot \sup (S_{2})
\end{cases}
\label{eq:multiplication}
 \end{equation*}

\noindent
and, as a particular case, we have 
 \begin{equation*}
\{k\}\boxdot S_2= S_2\boxdot\{k\}=\{ x \mid x = k \cdot s_{2}, s_{2} \in S_{2} \}
 \quad\text{with} \quad 
\begin{cases}
\inf (\{k\}\boxdot S_{2})= k\cdot\inf(S_2)\\
\sup(\{k\}\boxdot S_{2} )= k\cdot\sup(S_2)
\end{cases}
 \end{equation*}

{\it{Examples}}: 

$[0.1,0.6]\boxdot[0.8,0.9]=[0.08,0.54]$ because $0.1\cdot0.8=0.08$ and $0.6\cdot0.9=0.54$; 

$[0.1,0.6]\boxdot\{0.3\}=\{0.3\}\boxdot [0.1,0.6]=[0.03,0.18]$ because $0.3\cdot0.1=0.03$ and $0.3\cdot0.6=0.18$.\\

   \item
 {\bf{Division of sets}}

In our fusion context, the division of sets is not necessary since the DSm rules of combination (classic or hybrid ones) do not require a normalization procedure and thus a division operation. Actually, the  DSm rules require only addition and multiplication operations. We however give here the definition of division of sets only for the reader's  interest and curiosity. The division of sets is defined as follows:

If $0\notin S_2$, then
$ S_{1}\boxslash S_{2} \triangleq \{ x \mid x = s_{1}/s_{2}, 
s_{1} \in S_{1}, s_{2} \in S_{2}\}$  with 
 $\begin{cases}\inf(S_1\boxslash S_2)=\inf(S_1)/\sup(S_2)\\
 \sup(S_1\boxslash S_2)=\sup(S_1)/ \inf(S_2) \ \text{if} \ 0\not\in S_2\\
\sup(S_1\boxslash S_2)=+\infty \ \text{if}\ 0\in S_2
\end{cases}$

\noindent If $0\in S_2$, then  
$ S_{1}\boxslash S_{2}=[\inf(S_1)/\sup(S_2),+\infty)$

\noindent and as some particular cases, we have for $k\neq 0$, 
 \begin{equation*}
\{k\}\boxslash S_{2}=\{x \mid x = k/s_{2}, \text{where}\ s_{2} \in S_{2}\setminus \{0\}\}
 \quad\text{with} \quad 
 \begin{cases}
 \inf(\{k\}\boxslash S_2)=k/\sup(S_2)\\
 \sup(\{k\}\boxslash S_2)=k/ \inf(S_2)\end{cases}
 \label{eq:division2}
 \end{equation*}
 \noindent and if $0\in S_2$ then $\sup(\{k\}\boxslash S_{2})=+\infty$
 
One has also as some particular case for $k\neq 0$,
  \begin{equation*}
S_{2}\boxslash\{k\} =\{x \mid x = s_2/k, \text{where}\ s_{2} \in S_{2}\}
 \quad\text{with} \quad 
\begin{cases}\inf(S_{2}\boxslash\{k\})=\inf(S_2)/k\\
 \sup(S_{2}\boxslash\{k\})=\sup(S_2)/k\end{cases}
 \label{eq:division3}
 \end{equation*}

{\it{Examples}}: 

$[0.4,0.6]\boxslash[0.1,0.2]=[2,6]$ because $0.4/0.2=2$ and $0.6/0.1=6$; 

$[0.4,0.6]\boxslash\{0.4\}=[1,1.5]$ because $0.4/0.4=1$ and $0.6/0.4=1.5$; 

$\{0.8\}\boxslash[0.2,0.5]=[1.6,4]$ because $0.8/0.2=4$ and $0.8/0.5=1.6$; 

$[0,0.5]\boxslash[0.1,0.2]=[0,5]$: $[0,0.5]\boxslash\{0.4\}=[0,1.25]$ because $0/0.4=0$ and $0.5/0.4=1.25$; 

$[0.3,0.9]\boxslash[0,0.2]=[1.5,+\infty)$ because $0.3/0.2=1.5$ and $0\in (S_2=[0,0.2])$, $\sup([0.3,0.9]\boxslash[0,0.2])=+\infty$; 

$[0,0.9]\boxslash[0,0.2]=[0,+\infty)$: 

$\{0.7\}\boxslash[0,0.2]=[3.5,+\infty)$ because $0.7/0.2=3.5$ and $0\in (S_2=[0,0.2])$, $\sup(\{0.7\}\boxslash[0,0.2])=+\infty$; 

$\{0\}\boxslash[0,0.2]=[0,+\infty)$: $[0.3,0.9]\boxslash\{0\}=+\infty$: 

$[0,0.9]\boxslash\{0\}=+\infty$: 

$[0.2, 0.7]\boxslash[0, 0.8]= [0.25, +\infty)$.

\end{itemize}

These operations can be directly extended for any types of sets (not necessarily sub-unitary subsets as it will be showed in our general examples of section 6), but for simplicity, we will start the presentation in the following section only for sub-unitary subsets.\\

Due to the fact that the fusion of imprecise information must also be included in the unit interval $[0, 1]$ as it happens with the fusion of precise information, if the masses computed are less than $0$ one 
replaces them by $0$, and similarly if they are greater than $1$ one replaces them 
by $1$. For example (specifically in our fusion context): $[0.2, 0.4]\boxplus [0.5, 0.8]=[0.7, 1.2]$ will be forced to $[0.7, 1]$.

\section{Fusion of beliefs defined on single sub-unitary intervals }
%
\label{Sec:CIB}

\subsection{DSm rules of combination}

Let's now consider some given sources of information which are not able to provide us a specific/precise mass $m_{ij}\in [0,1]$, but only an interval centered\footnote{This interval centered assumption is not important actually but has been adopted here only for notation convenience.} in $m_{ij}$, i.e. $I_{ij}=[m_{ij}-\epsilon_{ij},m_{ij}+\epsilon_{ij}]$ where $0\leq\epsilon_{ij}\leq 1$ and $I_{ij}\subseteq [0,1]$ for all $1\leq i \leq k$ and $1\leq j\leq d$. The cases when $I_{ij}$ are half-closed or open are similarly treated.

\clearpage
\newpage
\noindent
{\bf{Lemma 1:}} if $A,B\subseteq [0,1]$ and $\alpha \in [0,1]$ then:

\begin{itemize}
\item[] 
$\begin{cases}
\inf (A \boxdot  B)=\inf (A)\cdot \inf (B)\\
\sup (A \boxdot  B)=\sup (A)\cdot \sup (B)
\end{cases}
\qquad\qquad
\begin{cases}
\inf (A \oplus B)=\inf (A) + \inf (B)\\
\sup (A \oplus B)=\sup (A) +  \sup (B)
\end{cases}$
\item[] $\begin{cases}
\inf (\alpha\cdot A)=\alpha\cdot \inf (A)\\
\sup (\alpha\cdot A)=\alpha\cdot \sup (A)
\end{cases}
\qquad\qquad\qquad\quad
\begin{cases}
\inf (\alpha + A)=\alpha +  \inf (A)\\
\sup (\alpha + A)=\alpha + \sup (A)
\end{cases}$
\end{itemize}

We can regard a scalar $\alpha $ as a particular interval $[\alpha,\alpha ]$,
 thus all operations of the previous lemma are reduced to multiplications and additions of sub-unitary intervals. Therefore, the DSm general rule \eqref{eq:DSmHkno1}, 
which operates (multiplies and adds) sub-unitary scalars, can be extended to operate sub-unitary intervals. The formula \eqref{eq:DSmHkno1} remains the same, but $m_i(X_i)$, $1\leq i\leq k$, are sub-unitary intervals $I_{ij}$. The mass matrix $\mathbf{M}$ is extended to:

\begin{equation*}
\inf(\mathbf{M})=
\begin{bmatrix}
m_{11}-\epsilon_{11} & m_{12}-\epsilon_{12} & \ldots & m_{1d}-\epsilon_{1d}\\
m_{21}-\epsilon_{21} & m_{22}-\epsilon_{22} & \ldots & m_{2d}-\epsilon_{2d}\\
 \ldots &  \ldots & \ldots & \ldots\\
m_{k1}-\epsilon_{k1} & m_{k2}-\epsilon_{k2} & \ldots & m_{kd}-\epsilon_{kd}\\
\end{bmatrix}
\end{equation*}
\begin{equation*}
\sup(\mathbf{M})=
\begin{bmatrix}
m_{11}+\epsilon_{11} & m_{12}+\epsilon_{12} & \ldots & m_{1d}+\epsilon_{1d}\\
m_{21}+\epsilon_{21} & m_{22}+\epsilon_{22} & \ldots & m_{2d}+\epsilon_{2d}\\
 \ldots &  \ldots & \ldots & \ldots\\
m_{k1}+\epsilon_{k1} & m_{k2}+\epsilon_{k2} & \ldots & m_{kd}+\epsilon_{kd}\\
\end{bmatrix}
\end{equation*}

\noindent
{\bf{Notations:}} Let's distinguish between DSm general rule for scalars, noted as usual $m_{\mathcal{M}(\Theta)}(A)$, or $m_i(X_i)$, etc., and the DSm general rule for intervals noted as 
$m_{\mathcal{M}(\Theta)}^{I}(A)$, or $m_i^{I}(X_i)$, etc. Hence, the DSm general rule for interval-valued masses is:

\begin{equation}
\inf(m_{\mathcal{M}(\Theta)}^I(A))\triangleq 
\phi(A)\Bigl[ S_1^{\inf} (A) + S_2^{\inf} (A) + S_3^{\inf}(A)\Bigr]
 \label{eq:DSmHkBisInf}
\end{equation}

\noindent
with
\begin{equation*}
S_1^{\inf}(A)\triangleq \sum_{\overset{X_1,X_2,\ldots,X_k\in D^\Theta}{(X_1\cap X_2\cap\ldots\cap X_k)=A}} \prod_{i=1}^{k} \inf(m_i^I(X_i))
\label{eq:S1inf}
\end{equation*}

 \begin{equation*}
S_2^{\inf}(A)\triangleq \sum_{\overset{X_1,X_2,\ldots,X_k\in\boldsymbol{\emptyset}}{ [\mathcal{U}=A]\vee [(\mathcal{U}\in\boldsymbol{\emptyset}) \wedge (A=I_t)]}} \prod_{i=1}^{k} \inf(m_i^I(X_i))
\label{eq:S2inf}
\end{equation*}
\begin{equation*}
S_3^{\inf}(A)\triangleq\sum_{\overset{X_1,X_2,\ldots,X_k\in D^\Theta}{\overset{(X_1\cup X_2\cup\ldots\cup X_k)=A}{\overset{(X_1\cap X_2\cap \ldots\cap X_k)\in\boldsymbol{\emptyset}}{}}}}  \prod_{i=1}^{k} \inf(m_i^I(X_i))
\label{eq:S3inf}
\end{equation*}
and
\begin{equation}
\sup(m_{\mathcal{M}(\Theta)}^I(A))\triangleq 
\phi(A)\Bigl[ S_1^{\sup} (A) + S_2^{\sup} (A) + S_3^{\sup}(A)\Bigr]
 \label{eq:DSmHkBisSup}
\end{equation}
\noindent
with
\begin{equation*}
S_1^{\sup}(A)\triangleq \sum_{\overset{X_1,X_2,\ldots,X_k\in D^\Theta}{(X_1\cap X_2\cap\ldots\cap X_k)=A}} \prod_{i=1}^{k} \sup(m_i^I(X_i))
\label{eq:S1sup}
\end{equation*}
 \begin{equation*}
S_2^{\sup}(A)\triangleq \sum_{\overset{X_1,X_2,\ldots,X_k\in\boldsymbol{\emptyset}}{ [\mathcal{U}=A]\vee [(\mathcal{U}\in\boldsymbol{\emptyset}) \wedge (A=I_t)]}} \prod_{i=1}^{k} \sup(m_i^I(X_i))
\label{eq:S2sup}
\end{equation*}
\begin{equation*}
S_3^{\sup}(A)\triangleq\sum_{\overset{X_1,X_2,\ldots,X_k\in D^\Theta}{\overset{(X_1\cup X_2\cup\ldots\cup X_k)=A}{\overset{(X_1\cap X_2\cap \ldots\cap X_k)\in\boldsymbol{\emptyset}}{}}}}  \prod_{i=1}^{k} \sup(m_i^I(X_i))
\label{eq:S3sup}
\end{equation*}

Actually formula \eqref{eq:DSmHkBisInf} results from applying the DSm hybrid rule for scalars to the matrix $\inf(\mathbf{M})$, while formula  \eqref{eq:DSmHkBisSup} results from applying the DSm hybrid rule for scalars to the matrix $\sup(\mathbf{M})$. The bounds of the DSm classic rule for  the free-DSm model are given for all $A\in D^\Theta$ by $S_1^{\inf}(A)$ and $S_1^{\sup}(A)$.
Combining  \eqref{eq:DSmHkBisInf} and \eqref{eq:DSmHkBisSup}, one gets directly:

\begin{equation}
m_{\mathcal{M}(\Theta)}^{I}(A)=[\inf m_{\mathcal{M}(\Theta)}^{I}(A), \sup m_{\mathcal{M}(\Theta)}^{I}(A)]
\end{equation}

Of course, the closeness of this interval to the left and/or to the right depends on the closeness of the combined intervals $I_{ij}$. If all of them are closed to the left, then $m_{\mathcal{M}(\Theta)}^{I}(A)$ is also closed to the left. But, if at least one is open to the left, then 
$m_{\mathcal{M}(\Theta)}^{I}(A)$ is open to the left. Similarly for the closeness to the right. Because one has $\forall i=1,\ldots,k$ and $\forall j=1,\ldots,d$ :
\begin{equation}
\lim_{\epsilon_{ij}\rightarrow 0}(\inf(\mathbf{M}))=\lim_{\epsilon_{ij}\rightarrow 0}(\sup(\mathbf{M}))= \mathbf{M}
\end{equation}

\noindent
It results the following theorem:\\

\noindent
{\bf{Theorem 1:}} $\forall A\in D^\Theta$, $\forall i=1,\ldots,k$ and $\forall j=1,\ldots,d$, one has:
\begin{equation}
\lim_{\epsilon_{ij}\rightarrow 0} m_{\mathcal{M}(\Theta)}^{I}(A) =[\lim_{\inf_{ij}} (A),\lim_{\sup_{ij}} (A)]
 \quad\text{with} \quad 
\begin{cases}
 \lim_{\inf_{ij}} (A) \triangleq \lim_{\epsilon_{ij}\rightarrow 0}(\inf (m_{\mathcal{M}(\Theta)}^{I}(A)))\\
 \lim_{\sup_{ij}} (A) \triangleq \lim_{\epsilon_{ij}\rightarrow 0}(\sup (m_{\mathcal{M}(\Theta)}^{I}(A)))
\end{cases}
\end{equation}

In other words, if all centered sub-unitary intervals converge to their corresponding mid points (the imprecision becomes zero), then the DSm rule for intervals converges towards the DSm rule for scalars.\\

Normally we must apply the DSm classical or  hybrid rules directly to the interval-valued masses, but this is equivalent to applying the DSm rules to the inferior and superior bounds of each mass.  If, after fusion, the sum of inferior masses is $< 1$ (which occurs all the time because combining incomplete masses one gets incomplete results) and the sum of superior masses is $\geq 1$ (which occurs all the time because combining paraconsistent masses one gets paraconsistent results), then there exist points in each resulted interval-valued mass such that their sum is 1 (according to a continuity theorem).

\subsection{Example with the DSm classic rule}

Let's take back the previous example (see section \ref{subsec:Examples2.2}), but let's now suppose the sources of information give at time $t$ imprecise generalized basic belief assignments, i.e. interval-valued masses centered in the scalars given in section \ref{subsec:Examples2.2}, of various radii according to table \ref{mytable4}.

\begin{table}[h]
\begin{equation*}
\begin{array}{|c|c|c|}
\hline
A\in D^\Theta & m_1^I(A) & m_2^I(A) \\
\hline
\theta_1 & [0.05,0.15] & [0.4,0.6]\\
\theta_2 & [0.1,0.3] & [0.1,0.5]\\
\theta_3 & [0.15,0.45] & [0,0.2]\\
\theta_1\cap\theta_2 & [0.2,0.6] & [0.05,0.15]\\
\hline
\end{array}
\end{equation*}
\caption{Inputs of the fusion with imprecise bba}
\label{mytable4}
\end{table}

Based on the free DSm model and the classical DSm rule applied to imprecise basic belief assignments following the method proposed in previous section, one has: 

$$m^I(\theta_1)=[0.05,0.15]\boxdot [0.4,0.6]=[0.020,0.090]$$
$$m^I(\theta_2)=[0.1,0.3]\boxdot  [0.1,0.5]=[0.010,0.150]$$
$$m^I(\theta_3)=[0.15,0.45] \boxdot [0,0.2]=[0,0.090]$$
$$m^I(\theta_1\cap \theta_3)=
[[0.05,0.15]\boxdot [0,0.2]]\boxplus
[[0.4,0.6]\boxdot [0.15,0.45]]=
[0,0.030]\boxplus [0.060,0.270]=[0.060,0.300]$$
$$m^I(\theta_2\cap \theta_3)=
[[0.1,0.3]\boxdot [0,0.2]]\boxplus
[[0.1,0.5]\boxdot [0.15,0.45]]=
[0,0.06]\boxplus [0.015,0.225]=[0.015,0.285]$$ 
\begin{align*}
m^I(\theta_1\cap\theta_2\cap \theta_3)&=
[[0.15,0.45]\boxdot [0.05,0.15]]\boxplus
[[0,0.2]\boxdot [0.2,0.6]]\\
&=
[0.0075,0.0675]\boxplus [0,0.12]\\
&=[0.0075,0.1875]
\end{align*}
\begin{align*}
m^I(\theta_1\cap \theta_2)& =
[[0.2,0.6]\boxdot [0.05,0.15]]\boxplus 
[[0.05,0.15]\boxdot  [0.05,0.15]]\boxplus
[[0.4,0.6]\boxdot [0.2,0.6]]\boxplus\\
& \qquad  
[[0.1,0.3]\boxdot [0.05,0.15]]\boxplus 
[[0.1,0.5]\boxdot [0.2,0.6]]\boxplus \\
& \qquad  
[[0.05,0.15]\boxdot  [0.1,0.5]]\boxplus 
[[0.4,0.6]\boxdot [0.1,0.3]]\\
&=
[0.010,0.90] \boxplus [0.0025,0.0225] \boxplus 
[0.08,0.36] \boxplus [0.005,0.045] \boxplus\\
& \qquad  
 [0.02,0.30] \boxplus [0.005,0.075] \boxplus [0.04,0.18] = [0.1625,1.0725]\equiv [0.1625,1]
\end{align*}
The last {\it{equality}} comes from the absorption of $[0.1625,1.0725]$ into $[0.1625,1]$ according to operations on sets defined in this fusion context. Thus, the final result of combination $m^I(.)=[m_1^I\oplus m_2^I](.)$ of these two imprecise sources of evidence is given in table \ref{mytable5}.
\begin{table}[h]
\begin{equation*}
\begin{array}{|c|c|}
\hline
A\in D^\Theta & m^I(A)=[m_1^I\oplus m_2^I](A) \\
\hline
\theta_1 & [0.020,0.090]\\
\theta_2 & [0.010,0.150]\\
\theta_3 & [0,0.090]\\
\theta_1\cap\theta_2 & [0.1625,1.0725\rightarrow 1]\\
\theta_1\cap\theta_3 & [0.060,0.300]\\
\theta_2\cap\theta_3 & [0.015,0.285]\\
\theta_1\cap\theta_2 \cap \theta_3& [0.0075,0.1875]\\
\hline
\end{array}
\end{equation*}
\caption{Fusion with DSm classic rule for free-DSm model}
\label{mytable5}
\end{table}

\noindent 
There exist some points, for example $0.03$, $0.10$. $0.07$, $0.4$, $0.1$, $0.2$, $0.1$ from 
the intervals $[0.020, 0.090]$, $\ldots$, $[0.0075, 0.1875]$ respectively such that 
their sum is 1 and therefore the admissibility of the fusion result holds. Note that this fusion process is equivalent to using the DSm classic rule for scalars for inferior limit and incomplete information (see table \ref{mytable6}), and the same rule for superior limit and paraconsistent information (see table \ref{mytable7}).
\begin{table}[h]
\begin{equation*}
\begin{array}{|c|c|c|c|}
\hline
A\in D^\Theta & m_1^{\inf}(A) & m_2^{\inf}(A) & m^{\inf}(A) \\
\hline
\theta_1 & 0.05 & 0.4 & 0.020\\
\theta_2 & 0.1 & 0.1 & 0.010\\
\theta_3 & 0.15 & 0 & 0\\
\theta_1\cap\theta_2 & 0.2 & 0.05 & 0.1625\\
\theta_1\cap\theta_3 & 0 & 0 & 0.060\\
\theta_2\cap\theta_3 & 0 & 0 & 0.015\\
\theta_1\cap\theta_2\cap\theta_3 & 0 & 0 & 0.0075\\
\hline
\end{array}
\end{equation*}
\caption{Fusion with DSm classic rule on lower bounds}
\label{mytable6}
\end{table}

\begin{table}[h]
\begin{equation*}
\begin{array}{|c|c|c|c|}
\hline
A\in D^\Theta & m_1^{\sup}(A) & m_2^{\sup}(A) & m^{\sup}(A) \\
\hline
\theta_1 & 0.15 & 0.6 & 0.090\\
\theta_2 & 0.3 & 0.5 & 0.150\\
\theta_3 & 0.45 & 0.2 & 0.090\\
\theta_1\cap\theta_2 & 0.6 & 0.15 & 1.0725\rightarrow 1\\
\theta_1\cap\theta_3 & 0 & 0 & 0.300\\
\theta_2\cap\theta_3 & 0 & 0 & 0.285\\
\theta_1\cap\theta_2\cap\theta_3 & 0 & 0 & 0.1875\\
\hline
\end{array}
\end{equation*}
\caption{Fusion with DSm classic rule on upper bounds}
\label{mytable7}
\end{table}

\subsection{Example with the DSm hybrid rule}

Then, assume at time $t+1$ one finds out for some reason that the free-DSm model has to be changed by introducing the constraint $\theta_1\cap\theta_2=\emptyset$ which involves also $\theta_1\cap\theta_2\cap\theta_3=\emptyset$. One directly applies the DSm hybrid rule for set to get the new belief masses:

\begin{align*}
m^I(\theta_1)&=[0.020,0.090]\boxplus [[0.05,0.15]\boxdot [0.05,0.15]]
\boxplus[ [0.4,0.6]\boxdot [0.2,0.6]]\\
&=[0.020,0.090]\boxplus [0.0025,0.0225] \boxplus [0.08,0.36]=[0.1025,0.4725]
\end{align*}

\begin{align*}
m^I(\theta_2)&=[0.010,0.150]\boxplus [[0.1,0.3]\boxdot[0.05,0.15]]
\boxplus [[0.1,0.5]\boxdot [0.2,0.6]]\\
&=
[0.010,0.150]\boxplus [0.005, 0.045]\boxplus [0.02,0.30] = [0.035, 0.495]
\end{align*}

\begin{align*}
m^I(\theta_3)&=[0,0.090]\boxplus [[0.15,0.45]\boxdot [0.05,0.15]]
\boxplus [[0,0.2]\boxdot [0.2,0.6]]\\
&=
[0,0.090]\boxplus [0.0075, 0.0675] \boxplus [0,0.12]=[0.0075,0.2775]
\end{align*}

\begin{align*}
m^I(\theta_1\cup\theta_2)&= [[02,0.6]\boxdot [0.05,0.15]]\boxplus
[[0.05,0.15]\boxdot [0.1,0.5]]\boxplus
[[0.4,0.6]\boxdot [0.1,0.3]]\\
&=
[0.010,0.090]\boxplus [0.005,0.075]\boxplus [0.04,0.18]=[0.055,0.345]
\end{align*}

$m^I(\theta_1\cap\theta_2)=m^I(\theta_1\cap\theta_2\cap\theta_3)=0$ by definition of empty masses (due to the choice of the hybrid model $\mathcal{M}$).
$m^I(\theta_1\cap\theta_3)=[0.060,0.300]$ and $m^I(\theta_2\cap\theta_3)=[0.015,0.285]$ remain the same. Finally, the result of the fusion of imprecise belief assignments for the chosen hybrid model $\mathcal{M}$, is summarized in table \ref{mytable8}.

\begin{table}[h]
\begin{equation*}
\begin{array}{|c|c|c|c|}
\hline
A\in D^\Theta & m^I(A)=[m^{\inf}(A),m^{\sup}(A)]\\
\hline
\theta_1 & [0.1025,0.4725]\\
\theta_2 & [0.035,0.495]\\
\theta_3 & [0.0075,0.2775]\\
\theta_1\cap\theta_2\overset{\mathcal{M}}{\equiv}\emptyset & [0,0]=0\\
\theta_1\cap\theta_3 & [0,060,0.300]\\
\theta_2\cap\theta_3 & [0.015,0.285]\\
\theta_1\cap\theta_2\cap\theta_3\overset{\mathcal{M}}{\equiv}\emptyset & [0,0]=0 \\
\theta_1\cup\theta_2 & [0.055,0.345] \\
\hline
\end{array}
\end{equation*}
\caption{Fusion with DSm hybrid rule for model $\mathcal{M}$}
\label{mytable8}
\end{table}

The admissibility of the fusion result still holds since there exist some points, for example $0.1$, $0.3$, $0.1$, $0$, $0.2$, $0.1$, $0, 0.2$ from the 
intervals $[0.1025, 0.4725]$, $\ldots$, $[0.055, 0.345] $ respectively such that 
their sum is 1. Actually in each of these examples there are infinitely many such 
groups of points in each respective interval whose sum is 1.  This can be 
generalized for any examples.

\section{Generalization of DSm rules for sets}
\label{Sec:GIB}

In this section, we extend the previous results on the fusion of admissible imprecise information defined only on single sub-unitary intervals to the general case where the imprecision is defined on sets. In other words, in previous section we dealt with admissible imprecise masses having the form $m^I(A)=[a,b] \subseteq [0,1]$, and now we deals with admissible imprecise masses having the form 
$m^I(A)=[a_1,b_1] \cup \ldots\cup [a_m,b_m] \cup (c_1,d_1)\cup  \ldots\cup (c_n,d_n)
\cup (e_1,f_1] \cup\ldots\cup (e_p,f_p] \cup [g_1,h_1)\cup \ldots\cup [g_q,h_q) \cup \{A_1,\ldots,A_r\}$ where all the bounds or elements involved into $m^I(A)$ belong to $[0,1]$. 

\subsection{General DSm rules for imprecise beliefs}

From our previous results, one can generalize the DSm classic rule from scalars to sets in the following way: $\forall A\neq\emptyset \in D^\Theta$,

\begin{equation}
m^I(A) = \underset{\underset{(X_1\cap X_2\cap\ldots\cap X_k)=A}{X_1,X_2,\ldots,X_k\in D^\Theta}}{\boxed{\sum}}
\underset{i=1,\ldots,k}{\boxed{\prod}} m_i^I(X_i)
\end{equation}

\noindent
where $\boxed{\sum}$ and $\boxed{\prod}$ represent the summation, and respectively product, of sets.\\

Similarly, one can generalize the DSm hybrid rule from scalars to sets in the following way:

\begin{equation}
m_{\mathcal{M}(\Theta)}^I(A)\triangleq 
\phi(A)\boxdot \Bigl[ S_1^I(A) \boxplus S_2^I(A) \boxplus S_3^I(A)\Bigr]
 \label{eq:DSmHkBisImprecise}
\end{equation}
\noindent
$\phi(A)$ is the {\it{characteristic emptiness function}} of the set $A$ and $S_1^I(A)$, $S_2^I(A)$ and $S_3^I(A)$ are defined by
\begin{equation}
S_1^I(A)\triangleq
\underset{\underset{(X_1\cap X_2\cap\ldots\cap X_k)=A}{X_1,X_2,\ldots,X_k\in D^\Theta}}{\boxed{\sum}}
\underset{i=1,\ldots,k}{\boxed{\prod}} m_i^I(X_i)
\label{eq:S1I}
\end{equation}
\begin{equation}
S_2^I(A)\triangleq 
\underset{\underset{[\mathcal{U}=A]\vee [(\mathcal{U}\in\boldsymbol{\emptyset}) \wedge (A=I_t)]}{X_1,X_2,\ldots,X_k\in\boldsymbol{\emptyset}}}{\boxed{\sum}}
\underset{i=1,\ldots,k}{\boxed{\prod}} m_i^I(X_i)
\label{eq:S2I}
\end{equation}

\begin{equation}
S_3^I(A)\triangleq
\underset{\underset{(X_1\cap X_2\cap\ldots\cap X_k)\in\boldsymbol{\emptyset} }{\underset{(X_1\cup X_2\cup\ldots\cup X_k)=A}{X_1,X_2,\ldots,X_k\in D^\Theta}}}{\boxed{\sum}}
\underset{i=1,\ldots,k}{\boxed{\prod}} m_i^I(X_i)
\label{eq:S3I}
\end{equation}

In the case when all sets are reduced to points (numbers), the set operations become normal operations with numbers; the sets operations are generalizations of numerical operations.

\subsection{Some lemmas and theorem}

\noindent
{\bf{Lemma 2:}} Let the scalars $a,b\geq 0$ and the intervals $I_1,I_2\subseteq [0,1]$, with $a\in I_1$ and $b\in I_2$. Then obviously $(a+b) \in I_1\boxplus I_2$ and $(a\cdot b) \in I_1\boxdot I_2$.\\

\noindent
Because in DSm rules of combining imprecise information, one uses only additions and subtractions of sets, according to this lemma if one takes at random a point of each mass set and one combines them using the DSm rules for scalars, the resulted point will belong to the resulted set from the fusion of mass sets using the DSm rules for sets.\\

\noindent
{\bf{Lemma 3:}} Let $\Theta=\{\theta_1,\theta_2,\ldots,\theta_n\}$ and $K\geq 2$ independent sources of information, and $d=\dim(D^\Theta)$. By combination of incomplete information in DSmT, one gets incomplete information.\\

\noindent
{\bf{Proof:}} Suppose the masses of the sources of information on $D^\Theta$ are for all $1\leq j \leq K$, represented by the mass-vector $\mathbf{m}_j=[m_{j_1},m_{j_2},\ldots,m_{j_d}]$ with $0\leq \sum_{r=1}^d m_{j_r} < 1$. According to the DSm network architecture, no matter what DSm rule of combination is applying (classic or hybrid), the sum of all resulted masses has the form:
\begin{equation}
\prod_{j=1}^K (m_{j_1} + m_{j_2} + \ldots + m_{j_d}) < (\underbrace{1\times 1\times\ldots\times1}_{K \;\text{times}})=1
\end{equation}

\noindent
{\bf{Lemma 4:}} By combination of paraconsistent information, one gets paraconsistent information.\\

\noindent
{\bf{Proof:}} Using the same notations and similar reasoning, one has for all $1\leq j\leq K$, $\mathbf{m}_j=[m_{j_1},m_{j_2},\ldots,m_{j_d}]$, with $\sum_{r=1}^d m_{j_r} > 1$. Then
$$\prod_{j=1}^K (m_{j_1} + m_{j_2} + \ldots + m_{j_d}) >  (\underbrace{1\times 1\times\ldots\times1}_{K \;\text{times}})=1$$

\noindent
{\bf{Lemma 5:}} Combining incomplete (sum of masses $< 1$) with complete (sum of masses $=1$) information, one gets incomplete information.\\

\noindent
{\bf{Lemma 6:}} Combining complete information, one gets complete information.\\

\noindent
{\bf{Remark:}} Combining incomplete with paraconsistent (sum of masses $> 1$) information can give any result. For example:
\begin{itemize}
\item If the sum of masses of the first source is 0.99 (incomplete) and the sum of masses of the second source is 1.01 (paraconsistent), then the sum of resulted masses is $0.99\times 1.01=0.9999$ (i.e. incomplete)
\item But if  the first is 0.9 (incomplete) and the second is 1.2 (paraconsistent), then the resulted sum of masses is $0.9\times1.2=1.08$ (i.e. paraconsistent).
\end{itemize}

We can also have: incomplete information fusionned with paraconsistent information and get complete information. For example: $0.8\times 1.25=1$.\\

\noindent
{\bf{Admissibility condition:}}  \\

An {\it{imprecise mass}} on $D^\Theta$ is considered {\it{admissible}} if there exist at least a point belonging to $[0,1]$ in each mass set such that the sum of these points is equal to 1 (i.e. complete information for at least a group of selected points).\\

\noindent Remark: A complete scalar information is admissible. Of course, for the incomplete scalar information and paraconsistent scalar information there can not be an admissibility condition, because by definitions the masses of these two types of informations do not add up to 1 (i.e. to the 
complete information).

\clearpage
\newpage
\noindent
{\bf{Theorem of Admissibility:}}  \\

Let  a frame $\Theta=\{\theta_1,\theta_2,\ldots,\theta_n\}$, with $n\geq 2$, its hyper-power set $D^\Theta$ with $\dim(D^\Theta)=d$, and $K\geq 2$ sources of information providing imprecise admissible masses on $D^\Theta$. Then, the resulted mass, after fusionning the imprecise masses of these sources of information with the DSm rules of combination, is also admissible.\\

\noindent
{\bf{Proof:}} Let $s_j$, $1\leq j \leq K$, be an imprecise source of information, and its imprecise admissible mass $\mathbf{m}_j^I=[m_{j_1}^I,m_{j_2}^I,\ldots,m_{j_d}^I]$. We underline that all $m_{j_r}^I$, for $1\leq r\leq d$, are sets (not scalars); if there is a scalar $\alpha$, we treat it as a set $[\alpha,\alpha]$. Because $\mathbf{m}_j^I$ is admissible, there exist the points (scalars in $[0,1]$)  $m_{j_1}^s \in m_{j_1}^I$, $m_{j_2}^s \in m_{j_1}^2$,\ldots,$m_{j_d}^s \in m_{j_d}^I$ such that $\sum_{r=1}^d m_{j_1}^s=1$. This property occurs for all sources of information, thus there exist such points $m_{j_r}^s$ for any $1\leq j\leq K$ and any $1\leq r\leq d$. Now, if we fusion, as a particular case, the masses of only these points, using DSm classic or hybrid rules, and according to lemmas, based on DSm network architecture, one gets complete information (i.e. sum of masses equals to 1). See also Lemma 2.


\subsection{An example with multiple-interval masses}
We present here a more general example with multiple-interval masses. For simplicity, this example is a particular case when the theorem of admissibility is verified by few points, which happen to be just on the bounders. More general and complex examples (not reported here due to space limitations), can be given and verified as well. It is however an extreme example, because we tried to comprise all kinds of possibilities which may occur in the imprecise or very imprecise fusion. So, let's consider a fusion problem over $\Theta=\{\theta_1,\theta_2\}$, two independent sources of information with the following imprecise admissible belief assignments
\begin{table}[h]
\begin{equation*}
\begin{array}{|c|c|c|}
\hline
A\in D^\Theta & m_1^I(A) & m_2^I(A) \\
\hline
\theta_1 & [0.1,0.2] \cup \{0.3\} & [0.4,0.5]\\
\theta_2 &(0.4,0.6)\cup [0.7,0.8] &  [0,0.4]\cup \{0.5,0.6\}\\
\hline
\end{array}
\end{equation*}
\caption{Inputs of the fusion with imprecise bba}
\label{mytablex1}
\end{table}

\noindent
Using the DSm classic rule for sets, one gets
\begin{align*}
m^I(\theta_1)&=([0.1,0.2] \cup \{0.3\})\boxdot [0.4,0.5] \\
&= ([0.1,0.2] \boxdot [0.4,0.5])\cup (\{0.3\}\boxdot [0.4,0.5] )\\
&= [0.04,0.10] \cup [0.12,0.15]\\
m^I(\theta_2)& =((0.4,0.6)\cup [0.7,0.8] )\boxdot ([0,0.4]\cup \{0.5,0.6\})\\
&= ((0.4,0.6)\boxdot [0,0.4])\cup ((0.4,0.6)\boxdot  \{0.5,0.6\}) \cup ([0.7,0.8]\boxdot [0,0.4]) \cup ([0.7,0.8]\boxdot \{0.5,0.6\})\\
&= (0,0.24)\cup (0.20,0.30) \cup (0.24,0.36)\cup [0,0.32] \cup [0.35,0.40] \cup [0.42,0.48] \\
&= [0,0.40] \cup [0.42,0.48]\\
m^I(\theta_1\cap \theta_2)&=
[([0.1,0.2] \cup \{0.3\})\boxdot([0,0.4]\cup \{0.5,0.6\})] \boxplus [[0.4,0.5]\boxdot ((0.4,0.6)\cup [0.7,0.8]) ]\\
&=[ ([0.1,0.2]\boxdot [0,0.4]) \cup ([0.1,0.2]\boxdot \{0.5,0.6\}) \cup ( \{0.3\}\boxdot [0,0.4]) \cup ( \{0.3\}\boxdot  \{0.5,0.6\})] \\
& \quad \boxplus [ ([0.4,0.5]\boxdot (0.4,0.6)) \cup ([0.4,0.5]\boxdot [0.7,0.8] ) ] \\
& = [[0,0.08]\cup [0.05,0.10]\cup [0.06,0.12] \cup [0,0.12] \cup \{0.15,0.18\}] \boxplus [(0.16,0.30)\cup[0.28,0.40]]\\
&= [[0,0.12]\cup \{0.15,0.18\}]\boxplus (0.16,0.40] \\
&=(0.16,0.52] \cup (0.31,0.55] \cup (0.34,0.58]\\
&=(0.16,0.58]
\end{align*}
\noindent
Hence finally the fusion admissible result is given by:
\begin{table}[h]
\begin{equation*}
\begin{array}{|c|c|}
\hline
A\in D^\Theta & m^I(A)= [m_1^I \oplus m_2^I](A) \\
\hline
\theta_1 & [0.04,0.10] \cup [0.12,0.15]\\
\theta_2 & [0,0.40] \cup [0.42,0.48] \\
\theta_1\cap \theta_2  & (0.16,0.58]\\
\theta_1\cup \theta_2 & 0 \\
\hline
\end{array}
\end{equation*}
\caption{Fusion result with the DSm classic rule}
\label{mytablex2}
\end{table}

\clearpage
\newpage

\noindent
If one finds out that $\theta_1\cap \theta_2 \overset{\mathcal{M}}{\equiv}\emptyset$ (this is our hybrid model $\mathcal{M}$ one wants to deal with), then one uses the DSm hybrid rule for sets \eqref{eq:DSmHkBisImprecise}: $m_{\mathcal{M}}^I(\theta_1\cap \theta_2)=0$ and $m_{\mathcal{M}}^I(\theta_1\cup \theta_2)= (0.16,0.58]$, the others imprecise masses are not changed. In other words, one gets now with DSm hybrid rule applied to imprecise beliefs:

\begin{table}[h]
\begin{equation*}
\begin{array}{|c|c|}
\hline
A\in D^\Theta & m_{\mathcal{M}}^I(A)= [m_1^I \oplus m_2^I](A) \\
\hline
\theta_1 & [0.04,0.10] \cup [0.12,0.15]\\
\theta_2 & [0,0.40] \cup [0.42,0.48] \\
\theta_1\cap \theta_2\overset{\mathcal{M}}{\equiv}\emptyset  & 0 \\
\theta_1\cup \theta_2 & (0.16,0.58]\\
\hline
\end{array}
\end{equation*}
\caption{Fusion result with the DSm hybrid rule for $\mathcal{M}$ }
\label{mytablex3}
\end{table}

Let's check now the admissibility conditions and theorem. For the source 1, there exist the precise masses $(m_1(\theta_1)=0.3) \in ([0.1,0.2] \cup \{0.3\})$ and $(m_1(\theta_2)=0.7) \in ((0.4,0.6)\cup [0.7,0.8])$ such that $0.3+0.7=1$. For the source 2, there exist the precise masses $(m_1(\theta_1)=0.4) \in ([0.4,0.5])$ and $(m_2(\theta_2)=0.6) \in ([0,0.4]\cup \{0.5,0.6\})$ such that $0.4+0.6=1$. Therefore both sources associated with $m_1^I(.)$ and $m_2^I(.)$ are admissible imprecise sources of information.\\

 It can be easily checked that the DSm classic fusion of $m_1(.)$ and $m_2(.)$ yields the paradoxical basic belief assignment $m(\theta_1)=[m_1\oplus m_2](\theta_1)=0.12$, $m(\theta_2)=[m_1\oplus m_2](\theta_2)=0.42$ and $m(\theta_1\cap \theta_2)=[m_1\oplus m_2](\theta_1\cap \theta_2)=0.46$.
One sees that the admissibility theorem is satisfied since $(m(\theta_1)=0.12)\in (m^I(\theta_1)=[0.04,0.10] \cup [0.12,0.15])$, $(m(\theta_2)=0.42)\in (m^I(\theta_2)=[0,0.40] \cup [0.42,0.48])$ and $(m(\theta_1\cap \theta_2)=0.46)\in (m^I(\theta_1\cap \theta_2)=(0.16,0.58])$ such that $0.12+0.42+0.46=1$. Similarly if one finds out that $\theta_1\cap\theta_2=\emptyset$, then one uses the DSm hybrid rule and one gets: $m(\theta_1\cap\theta_2)=0$ and $m(\theta_1\cup\theta_2)=0.46$; the others remain unchanged. The admissibility theorem still holds.

\section{ Conclusion}

In this paper, we proposed from the DSmT framework, a new general approach to combine, imprecise, uncertain and possibly paradoxist sources of information to cover a wider class of fusion problems. This work was motivated by the fact that in most of practical and real fusion problems, the information is rarely known with infinite precision and the admissible belief assignment masses, for each element of the hyper-power set of the problem, have to be taken/chosen more reasonably as sub-unitary (or as a set of sub-unitary) intervals rather than a pure and simple scalar values. This is a generalization of previous available works proposed in literature (mainly IBS restricted to TBM framework). One showed that it is possible to fusion directly interval-valued masses using the DSm rules (classic or hybrid ones) and the operations on sets defined in this work.  Several illustrative and didactic examples have been given throughout this paper to show the application of this new approach.  The method developed here can also combine incomplete and paraconsistent imprecise, uncertain and paradoxical sources of information as well. This approach (although focused here only on the derivation of imprecise basic belief assignments) can be extended without difficulty to the derivation of imprecise belief and plausibility functions as well as to imprecise pignistic probabilities according to the generalized pignistic transformation presented in \cite{Dezert_2004Book}. This work allows the DSmT to cover a wider class of fusion problems.

\end{document}